# RESEARCH ANNOUNCEMENT

# A NEW SERIES OF DENSE GRAPHS OF HIGH GIRTH

F. LAZEBNIK, V. A. USTIMENKO, AND A. J. WOLDAR

ABSTRACT. Let $k \geq 1$ be an odd integer, $t = \lfloor \frac{k+2}{4} \rfloor$, and $q$ be a prime power. We construct a bipartite, $q$-regular, edge-transitive graph $CD(k, q)$ of order $v \leq 2q^{k-t+1}$ and girth $g \geq k+5$. If $e$ is the the number of edges of $CD(k, q)$, then $e = \Omega(v^{1+\frac{1}{k-t+1}})$. These graphs provide the best known asymptotic lower bound for the greatest number of edges in graphs of order $v$ and girth at least $g$, $g \geq 5$, $g \neq 11, 12$. For $g \geq 24$, this represents a slight improvement on bounds established by Margulis and Lubotzky, Phillips, Sarnak; for $5 \leq g \leq 23$, $g \neq 11, 12$, it improves on or ties existing bounds.

## 1. INTRODUCTION

The missing definitions of graph-theoretical concepts which appear in this paper can be found in [5]. All graphs we consider are simple, i.e. undirected, without loops and multiple edges. Let $V(G)$ and $E(G)$ denote the set of vertices and the set of edges of $G$, respectively. $|V(G)| = v$ is called the *order* of $G$, and $|E(G)| = e$ is called the *size* of $G$. If $G$ contains a cycle, then the *girth* of $G$, denoted by $g = g(G)$, is the length of a shortest cycle in $G$. Some examples of graphs with large girth which satisfy some additional conditions have been known to be hard to construct and have turned out to be useful in different problems in extremal graph theory, in studies of graphs with a high degree of symmetry, and in the design of communication networks. There are many references on each of these topics. Here we mention just a few main books and survey papers which also contain extensive bibliographies. For extremal graph theory, see [5, 20]; for graphs with a high degree of symmetry, see [7]; for communication networks, see [8].

In this paper we present a new infinite series of regular bipartite graphs with edge-transitive automorphism group and large girth. More precisely, for each odd

Received by the editors October 26, 1993.
1991 *Mathematics Subject Classification.* Primary 05C35, 05C38.
This research was partially supported by NSF grant DMS-9115473.
The third author was additionally supported by NSF grant DMS-9304580 while at the Institute for Advanced Study, Princeton, NJ 08540.

integer $k \geq 1$ and any prime power $q$, we construct a bipartite, $q$-regular, edge-transitive graph $CD(k, q)$ of order at most $2q^{k-\lfloor \frac{k+2}{4} \rfloor +1}$ and girth at least $k+5$. Below we explain why these graphs are of interest.

1. Let $\mathcal{F}$ be a family of graphs. By $ex(v, \mathcal{F})$ we denote the greatest number of edges in a graph on $v$ vertices which contains no subgraph isomorphic to a graph from $\mathcal{F}$. Let $C_n$ denote the cycle of length $n \geq 3$. It is known (see [5, 6, 9]) that all graphs of order $v$ with more than $90kv^{1+\frac{1}{k}}$ edges necessarily contain a $2k$-cycle. Therefore $ex(v, \{C_3, C_4, \ldots, C_{2k}\}) \leq 90kv^{1+\frac{1}{k}}$. For a lower bound we know that $ex(v, \{C_3, C_4, \ldots, C_n\}) = \Omega(v^{1+\frac{1}{n-1}})$. The latter result follows from a theorem proved implicitly by Erdős (see [20]), and the proof is nonconstructive. As is mentioned in [20], it is unlikely that this lower bound is sharp, and several constructions support this remark for arbitrary $n$. For the best lower bounds on $ex(v, \{C_3, C_4, \ldots, C_{2s+1}\})$, $1 \leq s \leq 10$, see [1, 12–14, 20, 25, 27, 28]. For $s \geq 11$ and an infinite sequence of values of $v$, the best asymptotic lower bound $ex(v, \{C_3, C_4, \ldots, C_{2s+1}\} = \Omega(v^{1+\frac{2}{3s+3}})$ is provided by the family of Ramanujan graphs (see below).

Graphs $CD(k, q)$ show that for an infinite sequence of values of $v$,

$$ex(v, \{C_3, C_4, \ldots, C_{2s+1}\}) = \Omega(v^{1+\frac{2}{3s-3+\epsilon}}),$$

where $\epsilon = 0$ if $s$ is odd and $\epsilon = 1$ if $s$ is even. To our knowledge, this is the best known asymptotic lower bound for all $s$, $s \geq 2$, $s \neq 5$. For $s = 5$ a better bound $\Omega(v^{1+1/5})$ is given by the regular generalized hexagon.

2. Let $\{G_i\}$, $i \geq 1$, be a family of graphs such that each $G_i$ is an $r$-regular graph of increasing order $v_i$ and girth $g_i$. Following Biggs [2], we say that $\{G_i\}$ is a family of graphs with *large girth* if

$$g_i \geq \gamma \log_{r-1}(v_i)$$

for some constant $\gamma$. It is well known (e.g., see [5]) that $\gamma \leq 2$, but no family has been found for which $\gamma = 2$. For many years the only significant results in this direction were the theorems of Erdős and Sachs and its improvements by Sauer, Walther, and others (see p. 107 in [5] for more details and references), who, using nonconstructive methods, proved the existence of infinite families with $\gamma = 1$. The first explicit examples of families with large girth were given by Margulis [17] with $\gamma \approx 0.44$ for some infinite families with arbitrary large valency and $\gamma \approx 0.83$ for an infinite family of graphs of valency 4. The constructions were Cayley graphs of $SL_2(Z_p)$ with respect to special sets of generators. Imrich [11] was able to improve the result for an arbitrary large valency, $\gamma \approx 0.48$, and to produce a family of cubic graphs (valency 3) with $\gamma \approx 0.96$. In [4] a family of geometrically defined cubic graphs, so-called sextet graphs, was introduced by Biggs and Hoare. They conjectured that these graphs have large girth. Weiss [26] proved the conjecture by showing that for the sextet graphs (or their double cover) $\gamma \geq 4/3$. Then, independently, Margulis (see [18] and the references therein) and Lubotzky, Phillips, and Sarnak [16, 19] came up with similar examples of graphs with $\gamma \geq 4/3$ and arbitrary large valency (they turned out to be so-called Ramanujan graphs). In [3], Biggs and Boshier showed that $\gamma$ is exactly $4/3$ for the graphs from [16]. These are Cayley graphs of the group $PGL_2(Z_q)$ with respect to a set of $p+1$ generators, where $p$, $q$ are distinct primes congruent to 1 mod 4 with the Legendre symbol

$(\frac{p}{q}) = -1$.

In [13], Lazebnik and Ustimenko constructed the family of graphs $D(k, q)$ which give explicit examples of graphs with arbitrary large valency and $\gamma \geq \log_q(q-1)$. Their definition (see Section 2) and analysis are basically elementary. In [10] it was shown that for these graphs $\gamma = \log_q(q-1)$ for infinitely many values of $q$. Recently we discovered, with the aid of A. Schliep, that for $k \geq 6$, graphs $D(k, q)$ are disconnected, and, in fact, the number of connected components of these graphs grows exponentially with $k$. The main part of this paper is devoted to the analysis of these components. As they are all isomorphic for fixed $k$ and $q$, we denote any one of them by $CD(k, q)$. It will immediately follow that $\gamma \geq 4/3 \log_q(q-1)$ for the family of graphs $CD(k, q)$.

## 2. THE FAMILY $D(k, q)$

In this section we describe the graphs $D(k, q)$ and mention some of their properties. The interested reader is referred to [13] for additional information.

The construction of graphs $D(k, q)$ came out of an attempt to imitate, in some sense, the generalized $m$-gons which arise from rank two groups of Lie type. (In fact, when $k = 2$, 3, and 5, the graphs closely resemble the affine portions of a certain generalized triangle, quadrangle, and hexagon, respectively (see [12]).) The "coordinate" form in which the adjacency relations appear (see below) results from a technique of embedding Lie geometries in their corresponding Lie algebras and the notion of a blow-up of a graph (see [12, 21–24] for details).

Let $q$ be a prime power, and let $P$ and $L$ be two copies of the countably infinite dimensional vector space $V$ over $GF(q)$. Elements of $P$ will be called *points* and those of $L$ *lines*. In order to distinguish points from lines we introduce the use of parentheses and brackets: If $x \in V$, then $(x) \in P$ and $[x] \in L$. It will also be advantageous to adopt the notation for coordinates of points and lines introduced in [13]:

$$(p) = (p_1, p_{11}, p_{12}, p_{21}, p_{22}, p'_{22}, p_{23}, \ldots, p_{ii}, p'_{ii}, p_{i,i+1}, p_{i+1,i}, \ldots),$$

$$[l] = [l_1, l_{11}, l_{12}, l_{21}, l_{22}, l'_{22}, l_{23}, \ldots, l_{ii}, l'_{ii}, l_{i,i+1}, l_{i+1,i}, \ldots).$$

We now define an incidence structure $(P, L, I)$ as follows. We say point $(p)$ is incident to line $[l]$, and we write $(p)I[l]$ if the following relations on their coordinates hold:

(2.1)
$$\begin{aligned}
l_{11} - p_{11} &= l_1 p_1, \\
l_{12} - p_{12} &= l_{11} p_1, \\
l_{21} - p_{21} &= l_1 p_{11}, \\
l_{ii} - p_{ii} &= l_1 p_{i-1,i}, \\
l'_{ii} - p'_{ii} &= l_{i,i-1} p_1, \\
l_{i,i+1} - p_{i,i+1} &= l_{ii} p_1, \\
l_{i+1,i} - p_{i+1,i} &= l_1 p'_{ii}.
\end{aligned}$$

(The last four relations are defined for $i \geq 2$.) These incidence relations for $(P, L, I)$ become adjacency relations for a related bipartite graph. We speak now of the *incidence graph* of $(P, L, I)$, which has vertex set $P \cup L$ and edge set consisting of all pairs $\{(p), [l]\}$ for which $(p)I[l]$.

To facilitate notation in future results, it will be convenient for us to define
$p_{0,-1} = p_{-1,0} = l_{0,-1} = p_{1,0} = l_{0,1} = 0$, $p_{0,0} = l_{0,0} = -1$, $p'_{0,0} = l'_{0,0} = 1$,
$p_{0,1} = p_1$, $l_{1,0} = l_1$, $l'_{1,1} = l_{1,1}$, $p'_{1,1} = p_{1,1}$ and to rewrite (2.1) in the form:

(2.2)
$$\begin{aligned}
l_{ii} - p_{ii} &= l_1 p_{i-1,i}, \\
l'_{ii} - p'_{ii} &= l_{i,i-1} p_1, \\
l_{i,i+1} - p_{i,i+1} &= l_{ii} p_1, \\
l_{i+1,i} - p_{i+1,i} &= l_1 p'_{ii}, \\
\text{for} \quad i &= 0, 1, 2, \ldots.
\end{aligned}$$

Notice that for $i = 0$, the four conditions (2.2) are satisfied by every point and line, and for $i = 1$, the first two equations coincide and give $l_{1,1} - p_{1,1} = l_1 p_1$.

For each positive integer $k \geq 2$, we obtain an incidence structure $(P_k, L_k, I_k)$ as follows. First, $P_k$ and $L_k$ are obtained from $P$ and $L$, respectively, by simply projecting each vector onto its $k$ initial coordinates. Incidence $I_k$ is then defined by imposing the first $k-1$ incidence relations and ignoring all others. For fixed $q$, the incidence graph corresponding to the structure $(P_k, L_k, I_k)$ is denoted by $D(k, q)$. It is convenient to define $D(1, q)$ to be equal to $D(2, q)$. The properties of graphs $D(k, q)$ with which we are concerned are presented in the following

**Proposition 2.1.** *Let $q$ be a prime power and $k \geq 1$. Then*

(i) *$D(k,q)$ is a $q$-regular bipartite graph of order $2q^k$ ($2q^2$ for $k = 1$);*
(ii) *the automorphism group $\mathrm{Aut}(D(k, q))$ is transitive on points, lines, and edges;*
(iii) *for odd $k$, $g(D(k, q)) \geq k + 5$;*
(iv) *for odd $k$ and $q \equiv 1 \mod(\frac{k+5}{2})$, $g(D(k, q)) = k + 5$.* □

Proofs of parts (i), (ii), (iii) can be found in [13]; that of part (iv) in [10].

3. The family $CD(k, q)$

> *There is a crack in everything.*
> *That's how the light gets in.*
> — Leonard Cohen, *Anthem*

It turns out that for $k \geq 6$, graphs $D(k, q)$ are disconnected! Let $N_{k,q}$ be the number of connected components of $D(k, q)$.

**Lemma 3.1.** *Let $k \geq 1$, and let $t = \lfloor \frac{k+2}{4} \rfloor$. Then $N_{k,q} \geq q^{t-1}$.*

*Proof.* For $k = 1$ the statement is obvious, so we assume that $k \geq 2$. Let $u = (u_1, u_{11}, \ldots, u'_{tt}, \ldots)$ be a vertex of $D(k, q)$; it does not matter whether $u$ is a line or a point. For every $r$, $2 \leq r \leq t$, let

$$a_r = a_r(u) = \sum_{i=0}^{r}(u_{ii} u'_{r-i,r-i} - u_{i,i+1} u_{r-i,r-i-1}),$$

and $\vec{a} = \vec{a}(u) = (a_2, a_3, \ldots, a_t)$. Let $v$ be a vertex of $D(k, q)$ adjacent to $u$. We show that $\vec{a}(u) = \vec{a}(v)$. Without loss of generality, we may assume that $u$ is a point and $v$ is a line adjacent to $u$, say, $u = (p)$ and $v = [l]$. Then

$$a_r(p) = \sum_{i=0}^{r}(p_{ii} p'_{r-i,r-i} - p_{i,i+1} p_{r-i,r-i-1}).$$

Coordinates of $(p)$ can be expressed in terms of $p_1$ and the coordinates of $[l]$, use (2.2). When simplifying the resulting expression, we assume $l_{-1,-2} = l_{-1,-1} = l'_{-1,-1} = l_{-1,0} = 0$; the terms whose indices are out of range are multiplied by zeros, and therefore their appearance does not create a problem. Continuing, we eventually transform $a_r(p)$ into $a_r(l)$ (see [15] for the explicit computation).

Since $a_r(u) = a_r(v)$ for every $r$, $2 \leq r \leq t$, we get $\vec{a}(u) = \vec{a}(v)$ for any pair of adjacent vertices of $D(k,q)$. This implies that for any connected component $C$ of $D(k,q)$ and any vertices $x$, $y$ of $C$, $\vec{a}(x) = \vec{a}(y)$. Thus we may define $\vec{a}(C) = \vec{a}(v)$, where $v$ is a vertex of the connected component $C$.

Let us show that for every vector $\vec{c} = (c_2, c_3, \ldots, c_t) \in (GF(q))^{t-1}$ there exists a component $C$ of $D(k,q)$ such that $\vec{a}(C) = \vec{c}$. To do this, we just consider the following point $(p)$ in $D(k,q)$:

$$(p) = (0, 0, 0, 0, 0, p'_{22}, 0, 0, 0, p'_{33}, \ldots, 0, 0, 0, p'_{tt}, \ldots),$$

where $p'_{ii} = -c_i$ for all $i$, $2 \leq i \leq t$. Obviously, $\vec{a}(p) = \vec{c}$, and taking $C$ to be the connected component of $D(k,q)$ containing $(p)$, we obtain $\vec{a}(C) = \vec{c}$. Thus every $\vec{c} \in (GF(q))^{t-1}$ is "realizable" by a component of $D(k,q)$. Therefore $N_{k,q}$ is at least as large as $|(GF(q))^{t-1}| = q^{t-1}$. □

Due to the transitivity of $\mathrm{Aut}(D(k,q))$ on the the set of points of $D(k,q)$, (Proposition 2.1(ii)), all connected components of $D(k,q)$ are isomorphic graphs, and we denote any of them by $CD(k,q)$. We are ready to state the main result of this paper. Its proof is an immediate application of Proposition 2.1 and Lemma 3.1.

**Theorem 3.2.** *Let $k \geq 1$, $t = \lfloor \frac{k+2}{4} \rfloor$, $q$ be a prime power, $N_{k,q}$ be the number of connected components of $D(k,q)$, and $CD(k,q)$ be a connected component of $D(k,q)$. Then*

  (i) *$CD(k,q)$ is a bipartite, connected, $q$-regular graph of order $v = \frac{2q^k}{N_{k,q}} \leq 2q^{k-t+1}$;*
  (ii) *$\mathrm{Aut}(CD(k,q))$ acts transitively on points, lines, and edges of $CD(k,q)$;*
  (iii) *for odd $k$, the girth $g(CD(k,q)) \geq k+5$, and for $q \equiv 1 \mod(\frac{k+5}{2})$, $g(CD(k,q)) = k+5$.*
  (iv) *$e = \frac{vq}{2} = 2^{-1-\frac{1}{k}}(N_{k,q})^{\frac{1}{k}} v^{1+\frac{1}{k}}$, where $e$ is the size of $CD(k,q)$.* □

**Corollary 3.3.** *For $s \geq 2$, $ex(v, \{C_3, C_4, \ldots, C_{2s+1}\}) = \Omega(v^{1+\frac{2}{3s-3+\epsilon}})$, where $\epsilon = 0$ if $s$ is odd and $\epsilon = 1$ if $s$ is even.*

*Proof.* Set $2s = k+3$. Let $v$ and $e$ be the order and the size of the graph $CD(k,q)$. Then $v \leq 2q^{k-t+1}$, and $e = \frac{1}{2}vq \geq 2^{-1-\frac{1}{k-t+1}} v^{1+\frac{1}{k-t+1}}$. If $s$ is odd, then $k-t+1 = \frac{3s-3}{2}$; and if $s$ is even, then $k-t+1 = \frac{3s-2}{2}$. □

## Acknowledgment

We are indebted to Alexander Schliep, whose computer programs gave us an insight into the structure of the graphs $D(k,q)$.

## References

1. C. T. Benson, *Minimal regular graphs of girths eight and twelve*, Canad. J. Math. **18** (1966), 1091–1094.

2. N. L. Biggs, *Graphs with large girth*, Ars Combin. **25–C** (1988), 73–80.
3. N. L. Biggs and A. G. Boshier, *Note on the girth of Ramanujan graphs*, J. Combin. Theory Ser. B **49** (1990), 190–194.
4. N. L. Biggs and M. J. Hoare, *The sextet construction for cubic graphs*, Combinatorica **3** (1983), 153–165.
5. B. Bollobás, *Extremal graph theory*, Academic Press, London, 1978.
6. J. A. Bondy and M. Simonovits, *Cycles of even length in graphs*, J. Combin. Theory Ser. B **16** (1974), 97–105.
7. A. E. Brouwer, A. M. Cohen, and A. Neumaier, *Distance—regular graphs*, Springer-Verlag, Heidelberg and New York, 1989.
8. Fan R. K. Chung, *Constructing random-like graphs*, Probabilistic Combinatorics and its Applications, Proc. Sympos. Appl. Math., vol. 44, Amer. Math. Soc., Providence, RI, 1991.
9. R. J. Faudree and M. Simonovits, *On a class of degenerate extremal graph problems*, Combinatorica **3** (1983), 83–93.
10. Z. Füredi, F. Lazebnik, Á. Seress, V. A. Ustimenko, and A. J. Woldar, *Graphs of prescribed girth and bi-degree*, submitted.
11. W. Imrich, *Explicit construction of graphs without small cycles*, Combinatorica **2** (1984), 53–59.
12. F. Lazebnik and V. A. Ustimenko, *New examples of graphs without small cycles and of large size*, Europ. J. Combin. **14** (1993), 445–460.
13. F. Lazebnik and V. Ustimenko, *Explicit construction of graphs with an arbitrary large girth and of large size*, Discrete Appl. Math. (to appear).
14. F. Lazebnik, V. A. Ustimenko, and A. J. Woldar, *Properties of certain families of $2k$-cycle free graphs*, J. Combin. Theory Ser. B **60** (1994), 293–298.
15. \_\_\_\_\_\_, *A new series of dense graphs of large girth*, Rutcor Research Report RRR 99-93, December 1993.
16. A. Lubotzky, R. Phillips, and P. Sarnak, *Ramanujan graphs*, Combinatorica **8** (1988), 261–277.
17. G. A. Margulis, *Explicit construction of graphs without short cycles and low density codes*, Combinatorica **2** (1982), 71–78.
18. \_\_\_\_\_\_, *Explicit group-theoretical construction of combinatorial schemes and their application to the design of expanders and concentrators*, J. Problems of Inform. Trans. **24** (1988), 39–46; translation from Problemy Peredachi Informatsii **24** (January–March 1988), 51–60.
19. P. Sarnak, *Some applications of modular forms*, Cambridge Tracts in Math., vol. 99, Cambridge Univ. Press, Cambridge, 1990.
20. M. Simonovits, *Extremal graph theory*, Selected Topics in Graph Theory 2 (L. W. Beineke and R. J. Wilson, eds.), Academic Press, London, 1983, pp. 161–200.
21. V. A. Ustimenko, *Division algebras and Tits geometries*, Dokl. Akad. Nauk USSR **296** (1987), 1061–1065. (Russian)
22. \_\_\_\_\_\_, *A linear interpretation of the flag geometries of Chevalley groups*, Kiev Univ., Ukrain. Mat. Zh. **42** (March 1990), 383–387.
23. \_\_\_\_\_\_, *On the embeddings of some geometries and flag systems in Lie algebras and superalgebras*, Root Systems, Representation and Geometries, IM AN UkrSSR, Kiev, 1990, pp. 3–16.
24. \_\_\_\_\_\_, *On some properties of geometries of the Chevalley groups and their generalizations*, Investigation in Algebraic Theory of Combinatorial Objects (I. A. Faradzev, A. A. Ivanov, M. H. Klin, and A. J. Woldar, eds.) Kluwer, Dordrecht, 1991, pp. 112–121.
25. V. A. Ustimenko and A. J. Woldar, *An improvement on the Erdős bound for graphs of girth* 16, Proceedings of International Conference in Algebra, Barnaul, Russia, September 1991, Contemp. Math., Amer. Math. Soc., Providence, RI (to appear).

26. A. I. Weiss, *Girth of bipartite sextet graphs*, Combinatorica **4** (1984), 241–245.
27. R. Wenger, *Extremal graphs with no $C^4$, $C^6$, or $C^{10}$'s*, J. Combin. Theory Ser. B **52** (1991), 113–116.
28. A. J. Woldar and V. A. Ustimenko, *An application of group theory to extremal graph theory*, Group Theory, Proceedings of the Ohio State-Denison Conference, World Scientific, Singapore, 1993.

Department of Mathematical Sciences, University of Delaware, Newark, Delaware 19716
*E-mail address*: `fellaz@math.udel.edu`

Department of Mathematics and Mechanics, University of Kiev, Kiev 252127, Ukraine
*E-mail address*: `vau@rpd.univ.kiev.ua`

Department of Mathematical Sciences, Villanova University, Villanova, Pennsylvania 19085
*E-mail address*: `woldar@vill.edu`